\documentstyle[amscd,amssymb,verbatim,diagrams,12pt]{amsart}
\newcommand{\Si}{\Sigma}

\newcommand{\ZZ}{{\cal Z}}

\newcommand{\Cone}{\operatorname{Cone}}

\newcommand{\OO}{{\cal O}}

\newcommand{\BB}{{\cal B}}

\newcommand{\G}{{\Bbb G}}

\newcommand{\lan}{\langle}
\newcommand{\ran}{\rangle}
\newcommand{\Coh}{\operatorname{Coh}}
\newcommand{\GG}{{\cal G}}
\newcommand{\CC}{{\cal C}}

\newcommand{\Spec}{\operatorname{Spec}}

\newcommand{\ga}{\gamma}
\newcommand{\de}{\delta}

\renewcommand{\ker}{\operatorname{ker}}

\newcommand{\qcoh}{\operatorname{qcoh}}
\numberwithin{equation}{subsection}

\newcommand{\GL}{\operatorname{GL}}

\newtheorem{thm}{Theorem}[subsection]

\newtheorem{lem}[thm]{Lemma}

{  \theoremstyle{definition}

\newtheorem{ex}[thm]{Example}

}

\newcommand{\Pf}{\noindent {\it Proof}}
\newcommand{\id}{\operatorname{id}}

\newcommand{\ov}{\overline}

\newcommand{\rk}{\operatorname{rk}}
\newcommand{\rig}{\operatorname{rig}}

\newcommand{\FF}{{\cal F}}

\newcommand{\XX}{{\cal X}}

\newcommand{\PP}{{\cal P}}

\newcommand{\VV}{{\cal V}}
\newcommand{\SS}{{\cal S}}

\newcommand{\Ext}{\operatorname{Ext}}

\newcommand{\Res}{\operatorname{Res}}

\newcommand{\Aut}{\operatorname{Aut}}

\renewcommand{\a}{\alpha}
\renewcommand{\b}{\beta}
\newcommand{\om}{\omega}
\newcommand{\De}{\Delta}
\newcommand{\la}{\lambda}

\newcommand{\C}{{\Bbb C}}

\newcommand{\Z}{{\Bbb Z}}

\newcommand{\Ga}{\Gamma}

\newcommand{\wt}{\widetilde}
\newcommand{\ot}{\otimes}

\newcommand{\sub}{\subset}
\newcommand{\ed}{\qed\vspace{3mm}}

\newcommand{\Qcoh}{\operatorname{Qcoh}}

\newcommand{\tot}{\operatorname{tot}}
\newcommand{\MF}{\operatorname{MF}}

\title{Fundamental matrix factorization in the FJRW-theory revisited}
\author{Alexander Polishchuk}
\thanks{Supported in part by the NSF grant DMS-1700642 and by the Russian Academic Excellence Project `5-100'.}

\begin{document}

\maketitle
\begin{abstract}
We present an improved construction of the fundamental matrix factorization in the FJRW-theory
given in \cite{PV}. The revised construction
is coordinate-free and works for a possibly nonabelian finite group of symmetries.
One of the new ingrediants is the category of dg-matrix factorizations over a dg-scheme.
\end{abstract}

\section*{Introduction}

This short note is supposed to clarify the construction of the cohomological field theory associated
with a quasihomogeneous polynomial $W$ and its finite group of symmetries $G$. Such a cohomological field theory,
called the {\it FJRW-theory} was first proposed in \cite{FJR}. Then, in \cite{PV} a different construction, based
on categories of matrix factorizations, was given (conjecturally, the two constructions give the same cohomological field
theory). 

The approach of \cite{PV} is based on constructing certain {\it fundamental matrix factorizations}
which live over the product of certain finite coverings of $\ov{M}_{g,n}$ (the moduli of $\Ga$-spin structures) with
affine spaces. It is this construction that we aim to clarify.  
More precisely, we would like to present the construction in such a way that it would be analogous
to the construction of Ciocan-Fontanine and Kapranov of the virtual fundamental class in Gromov-Witten
theory via dg-manifolds (see \cite{CF-K}). The second goal that we achieve is to present the construction without
using coordinates on the vector space $V$ on which $W$ lives. This has an additional bonus that we can handle the
case when the group $G$ is not necessarily commutative (but still finite). 


The construction of \cite{PV} of the fundamental matrix factorization over $\SS\times \prod_i V^{\ga_i}$, where $\SS$ is the moduli
space of (rigidified) $\Ga$-spin structures with some markings (see Sec.\ \ref{setup-sec} for details)
roughly has the following two steps. In Step 1 one considers the object $R\pi_*(\VV)$ in the derived category
$D(\SS)$, where $\pi:\CC\to \SS$ is the universal curve, $\VV$ is the underlying vector bundle of the universal
$\Ga$-spin structure, and then equips it with some additional structure. In Step 2 one realizes $R\pi_*(\VV)$ by a $2$-term
complex $[A\to B]$, such that there is a morphism 
$$Z:X=\tot(A)\to \prod_i V^{\ga_i}$$
and a Koszul matrix factorization of $Z^*(\sum W_i)$, where $W_i=W|_{V^{\ga_i}}$. 
Then the fundamental matrix factorization is obtained by taking
its push-forward with respect to the morphism $(p,Z):X\to \SS\times V^{\ga_i}$, where $p:X\to \SS$ is the projection.
Note that here the space $X$ is non-canonical, so one has to check independence on the choices made.

The main idea of the present paper is to change the conceptual framework slightly by observing that in fact one gets a 
{\it dg-matrix factorization} on a {\it dg-scheme} over 
$\SS\times \prod_i V^{\ga_i}$ (the terminology is explained in Sec.\ \ref{mf-sec}).
Namely, for a non-negatively graded complex of vector bundles $C^\bullet$ over $\SS$, one can define the corresponding
dg-scheme over $\SS$,
$$[C^\bullet]:=\Spec(S^\bullet(C^\bullet)^\vee).$$
In our case we consider the dg-scheme 
$$\XX:=[R\pi_*(\VV)].$$ 
More concretely, if we realize $\VV$ by a $2$-term complex
$\VV=[A\to B]$ then our dg-scheme is realized by the sheaf of dg-algebras 
$$\OO_{\XX,[A\to B]}:=S^\bullet(B^\vee\to A^\vee).$$
Then we interpret the additional structure on $R\pi_*(\VV)$ coming from the universal $\Ga$-spin structure as
a structure of a dg-matrix factorization on the structure sheaf of $\XX$. More precisely, we get a morphism
$$Z_{\XX}:\XX\to \prod_i V^{\ga_i}$$
and a function of degree $-1$, $f_{-1}\in \OO_{\XX,[A\to B]}^{-1}$, such that
$$d(f_{-1})=-Z_{\XX}^*(\sum W_i).$$
Now 
the fundamental matrix factorization is obtained as the push-forward of $(\OO_\XX, d+f_{-1}\cdot\id)$
with respect to the morphism $\XX\to \SS\times \prod_i V^{\ga_i}$.

The connection with the original approach is the following: for each presentation $\VV=[A\to B]$, for which the
first construction works, there is a morphism $q:\XX\to X=\tot(A)$, such that $Z\circ q=Z_{\XX}$, 
and an isomorphism of the push-forward
$q_*(\OO_\XX,d+f_{-1}\cdot\id)$ with the Koszul matrix factorization of $Z^*(\sum W_i)$ constructed through the first approach.

The second technical improvement we present is in the construction of $f_{-1}$. The idea is to work systematically 
with the categories of sheaves over pairs (scheme, closed subscheme)
to deal with non-functoriality of the cone construction (such categories fit into the framework of
Lunts's poset schemes in \cite{lunts-poset}) .
Namely, we work with the enhancement of 
the usual push-forward with respect to the projection $\pi:\CC\to \SS$ to a morphism
of pairs $(\CC,\Sigma)\to (\SS,\SS)$, where $\Sigma\sub \CC$ is the union of the images of the universal marked points
(see Sec.\ \ref{pairs-sec}).

Throughout this work the ground field is $\C$.

\noindent
{\it Acknowledgments}. I am grateful to Felix Janda and Yongbin Ruan for organizing the RTG Conference on Witten's
$r$-spin class and related topics in January 2017, where the results of this note were first presented.
I also thank Institut Mathematique Jussieu and Institut des Hautes Etudes Scientifiques for hospitality and excellent
working conditions during preparation of this paper.

\section{Matrix factorizations over dg-schemes}\label{mf-sec}

\subsection{Definition}

We consider dg-schemes in the spirit of \cite{CF-K}. We fix a space $S$ (a scheme or a stack),
and consider the structure sheaf of a dg-scheme over $S$ to be a sheaf 
$(\OO^\bullet_X,d)$ of $\Z_-$-graded commutative dg-algebras over $\OO_S$
(one can make a restriction $\OO^0_X=\OO_S$, but it is not really necessary).

Given a function $f_0\in \OO^0_X$ we can consider the category of (quasicoherent) {\it dg-matrix factorizations} of $f_0$.
By definition, these are $\Z/2$-graded complexes of sheaves $P=P^{\ov{0}}\oplus P^{\ov{1}}$ together with
a (quasicoherent) $\OO_X^\sharp$-module structure, such that $\OO_X^i\cdot P^{\ov{a}}\sub P^{\ov{i+a}}$.
In addition $P$ is equipped with an odd differential $\de$
satisfying the Leibnitz identity
$$\de(\phi\cdot p)=d(\phi)\cdot p+(-1)^k\phi\de(p),$$
for $\phi\in \OO_X^k$, $p\in P$, 
and the equation $\de^2=f_0\cdot\id_P$.

\begin{ex}\label{cob-ex} Given an element $f_{-1}\in \OO^{-1}_X$, such that $d(f_{-1})=f_0$, 
we get a structure of a dg-matrix factorization on $\OO^\bullet_X$ by setting 
$$\de(\phi)=d(\phi)+f_{-1}\cdot\phi.$$
(In checking that $\de^2=0$ one has to use the fact that $f_{-1}^2=0$.)
\end{ex}


The above example can be obtained from the following more general operation.
Suppose we are given a function $f_0\in \OO^0_X$ and a dg-matrix factorization $(P,\de)$ of $f_0$.
Then for any $f_{-1}\in \OO^0_X$ we can change the differential $\de$ to $\de+f_{-1}\cdot \id_P$.
Then $(P,\de+f_{-1}\cdot \id_P)$ will be a dg-matrix factorization of $f_0+d(f_{-1})$.


\subsection{Positselski's framework of quasicoherent CDG-algebras}

More generally, we can assume that $f_0$ a section in $\OO^0_X\otimes L$, where $L$ is a locally free $\OO^0_X$-module
of rank $1$. The theory of the corresponding categories of dg-matrix factorizations fits into the framework of
quasicoherent CDG-algebras developed by Positselski (see \cite[Sec.\ 1]{EP}).

With the data $(\OO^\bullet_X, L, f_0)$ as above we can associate a quasicoherent CDG-algebra
$$\BB:=\bigoplus_{n\in \Z} \OO^\bullet_X\otimes_{\OO^0_X} L^{\otimes n}[-2n],$$
with the natural structure of a complex of sheaves (i.e., the $\Z$-grading and the differential $d$),
the natural product and the global curvature element given by $f_0\in \OO^0\ot L\sub \BB_2$.

Now a quasicoherent dg-matrix factorization is a quasicoherent DG-module over $\BB$, i.e.,
a graded $\BB$-module $M=\bigoplus_n M_n$, equipped with a differential $\de=\de_M$ such that $\de^2=f_0\cdot \id_M$
and $\de$ satisfies the Leibnitz identity with respect to the $\BB$-action.
Note that such a DG-module necessarily has 
$$M_{n+2}\simeq M_n\otimes L,$$
so it is determined by the components $M_0$ and $M_1$, and we get the structure of a dg-matrix factorization
on $M_0\oplus M_1$.

There are several exotic derived categories associated to a quasicoherent CDG-algebra.
The one that is most relevant for the theory of dg-matrix factorizations is the category
$$\qcoh-\MF_{ffd}(f_0):=
D^{co}(\BB-\qcoh_{ffd})\simeq D^{co}(\BB-\qcoh_{fl})\simeq D^{abs}(\BB-\qcoh_{fl}),$$
where the superscripts "abs" and "co" refer to "absolute" and "coderived", while the subscripts "fl" and "ffd"
mean "flat" and "finite flat dimension" (see \cite[Sec.\ 1]{EP}). 





Assume that $f:(X,\OO^{\bullet}_X)\to (Y,\OO^\bullet_Y)$ is a morphism of finite flat dimension, 
$L$ is a locally free $\OO^0_Y$-module of rank $1$, $W_0$ is a section of $L$.
Then we have the induced section $f^*W_0$ of $f^*L$.
In this situation we have the push-forward functor (see \cite[Prop.\ 1.9]{EP}) 
$$Rf_*:\qcoh-\MF_{ffd}(f^*W_0)\to \qcoh-\MF_{ffd}(W_0).$$

\subsection{Koszul matrix factorizations as push-forwards}
 
Let $V$ be a vector bundle over $X$, and suppose we have
sections $\a\in H^0(X,V^\vee)$, $\b\in H^0(X,V)$. With these data one associates a Koszul matrix factorization 
$\{\a,\b\}$ of $W=\lan \a,\b\ran$, whose underlying super-vector bundle is $\bigwedge^\bullet(V)$.
On the other hand, we have the derived zero locus of $\b$, $\ZZ(\b)\to X$, which corresponds to the dg-algebra
given by the Koszul complex of $\b$:
$$\OO_{\ZZ(\b)}=({\bigwedge}^\bullet(V),d=\iota_\b).$$
Now we can view $\a$ as a function of degree $-1$ on $\ZZ(\b)$ such that $d(\a)$ is the pull-back of $W$.
Thus, by definition, $\{\a,\b\}$ is the push-forward of the dg-matrix factorization $(\OO_{\ZZ(\b)},d+\a\cdot\id)$.

This explains why in the case when $\b$ is a regular section of $V$, the Koszul matrix factorization $\{\a,\b\}$ is equivalent
to the push-forward of the structure sheaf on the usual zero locus of $\b$.




\section{Trace maps via morphisms of pairs}\label{pairs-sec}

\subsection{Sheaves on pairs}

Let $\iota:Y\to X$ be a closed embedding.

We consider a very simple poset scheme in the sense of \cite{lunts-poset} for the poset consisting of two elements $\a>\b$,
so that $X_\a=Y$ and $X_\b=X$. Then a quasicoherent sheaf on this poset scheme is a triple
$(\FF_\a,\FF_\b,\phi)$, with $\FF_\a\in\Qcoh(Y)$, $\FF_\b\in\Qcoh(X)$ and 
$\phi:\FF_\b\to \iota_*\FF_\a$ is a morphism. We denote by $\Qcoh(X,Y)$ this abelian category,
and by $\Coh(X,Y)$ its subcategory corresponding to $\FF_\a\in\Coh(Y)$, $\FF_\b\in\Coh(X)$.

Note that the derived category $D^b\Coh(X,Y)$ has a natural monoidal structure given by the tensor product, so we can also define symmetric powers of objects in $D^b\Coh(X,Y)$.

Given a morphism of pairs $f:(X,Y)\to (X',Y')$ we have a natural derived push-forward morphism
$$Rf_*:D^+\Qcoh(X,Y)\to D^+\Qcoh(X',Y').$$

The push-forward is compatible with the tensor products in the usual way: we have natural morphisms
\begin{equation}\label{push-forward-tensor-eq}
Rf_*(F)\ot Rf_*(G)\to Rf_*(F\ot G), \ \ S^\bullet Rf_*(F)\to Rf_* S^\bullet(F).
\end{equation}

We have a fully faithful exact embedding $j_!:D\Qcoh(X)\to D\Qcoh(X,Y)$ sending $\GG$ to
$\FF_\a=0$, $\FF_\b=\GG$. There is a right adjoint functor to it (see \cite{lunts-poset}),
$$Rj^!:D^+\Qcoh(X,Y)\to D^+\Qcoh(X),$$
which is defined as the right derived functor of the functor 
$$j^!:\Qcoh(X,Y)\to \Qcoh(X): \FF_\bullet\mapsto \ker(\FF_\b\to \iota_*\FF_\a).$$
Note that objects $\FF_\bullet\in\Qcoh(X,Y)$, such that $\FF_\b\to \iota_*\FF_\a$ is surjective,
are acyclic with respect to $j^!$. Furthermore, 
every object of $\Qcoh(X,Y)$ has a canonical resolutions by such acyclic objects: 
$$0\to (\FF_\a,\FF_\b)\to (\FF_\a, \FF_\b\oplus\iota_*\FF_\a)\to (0, \iota_*\FF_\a)\to 0$$
Computing $Rj^!$ using these resolutions has a very simple interpretation: given a complex
$(\FF^\bullet_\a,\FF^\bullet_\b)$ over $\Qcoh(X)$, the functor $Rj^!$ sends it to the complex
$$\Cone(\FF^\bullet_\b\to \iota_*\FF^\bullet_\a)[-1].$$
In particular, there is a natural exact triangle
$$Rj^!(\FF^\bullet_\a,\FF^\bullet_\b)\to \FF^\bullet_\b\to \iota_*\FF\bullet_\a\to\ldots$$



We also have the following compatibility between $Rj^!$ and the push-forward.

\begin{lem}\label{i!-comm-lem}
Let $f:(X,Y)\to (X',Y')$ be a morphism of pairs. Assume that there exists a finite open covering of $X$, affine over $X'$.
Then for $\FF\in D^+\Qcoh(X,Y)$
we have a natural isomorphism
\begin{equation}\label{i!-f*-eq}
Rj^!Rf_*(\FF)\simeq Rf_*Rj^!(\FF)
\end{equation}
in $D^+\Qcoh(X')$.
\end{lem}

\Pf . Let us choose a quasi-isomorphism $\FF\to \wt{\FF}$, such that all $\wt{\FF}^i_\a$ and $\wt{\FF}^i_\b$ are 
$f_*$-acyclic (this can be done using Cech resolutions).
Then the left-hand side of \eqref{i!-f*-eq} is represented by the complex
$$\Cone(f_*\wt{\FF}_\b\to \iota_*f_*\wt{\FF}_\a)[-1].$$
On the other hand, the terms of $\Cone(\wt{\FF}_\b\to \iota_*\wt{\FF}_\a)[-1]$ are also $f_*$-acyclic, so
the right-hand side of \eqref{i!-f*-eq} is represented by the complex
$$f_*\Cone(\wt{\FF}_\b\to \iota_*\wt{\FF}_\a)[-1],$$
which is isomorphic to the one above.
\ed

\subsection{Differentials on curves}

Let $\pi:\CC\to\SS$ be a family of stable curves, $p_i:\SS\to \CC$, $i=1,\ldots,r$, be sections of $\pi$,
such that $\pi$ is smooth along their images, and let $\Si=\sqcup_i p_i(\SS)$. We view $(\CC,\Si)$ as a poset
scheme and consider the corresponding category $\Coh(\CC,\Si)$ whose objects are collections $(F,(F_i),(f_i))$,
where $F$ is a coherent sheaf on $\CC$, $F_i$ is a coherent sheaf on $\SS$ and $f_i:F\to p_{i*}F_i$ is a morphism.
Sometimes we will omit the morphisms $(f_i)$ from the notation and just write $(F,(F_i))$.

Set $\om_{\CC/\SS}^{\log}=\om_{\CC/\SS}(\Si)$.
Recall that we have natural residue maps 
$$\Res_{\Si}: \om_{\CC/\SS}^{\log}|_{\Si}\rTo{\sim} \OO_{\Si},$$
so that $\ker(\Res_{\Si})$ is identified with $\om_{\CC/\SS}$.
Thus, we can view the triple 
$$[\om_{\CC/\SS}^{\log},\Si]:=(\om_{\CC/\SS}^{{\log}},\OO_{\Si},\Res_{\Si})$$ 
as an object of the category
$\Coh(\CC,\Si)$. Furthermore, we have 
$$Rj^![\om_{\CC/\SS}^{\log},\Si]\simeq \om_{\CC/\SS}.$$
Note that we have a morphism of pairs
\begin{equation}\label{projection-poset-eq}
\pi:(\CC,\Si)\to (\SS,\SS).
\end{equation}
By Lemma \ref{i!-comm-lem}, the object $R\pi_*[\om_{\CC/\SS}^{\log},\Si]$ satisfies
\begin{equation}\label{j!-om-log-eq}
Rj^!R\pi_*[\om_{\CC/\SS}^{\log},\Si]\simeq R\pi_*\om_{\CC/\SS}.
\end{equation}

Note also that we have a morphism of exact triangles (which will be used later)
\begin{equation}\label{diff-res-diagram}
\begin{diagram}
\oplus_{i=1}^r\OO_\SS[-1]&\rTo{}&R\pi_*(\om_{\CC/\SS})&\rTo{}&R\pi_*(\om^{\log}_{\CC/\SS})&\rTo{}& \oplus_{i=1}^r\OO_\SS\\
\dTo{}&&\dTo{}&&\dTo{}&&\dTo{t}\\
\OO_S[-1]&\rTo{\id}&\OO_S[-1]&\rTo{}&0 &\rTo{}& \OO_\SS
\end{diagram}
\end{equation}
where $t$ is given by the summation.

The above constructions also work in the case of a family of orbicurves with stable coarse moduli spaces.

\section{Fundamental matrix factorization}

\subsection{Setup and the moduli spaces of $\Ga$-spin structures}\label{setup-sec}

Let us recall the setup of the FJRW theory (see \cite{FJR}, \cite{PV}), or rather its slight generalization to noncommutative
finite groups of symmetries (as in \cite{FJR2}).

We start with a finite-dimensional vector space $V$ equipped with an effective $\G_m$-action called the {\it $R$-charge},
such that all the weights of this action on $V$ are positive.
We denote the corresponding subgroup in $\GL(V)$ by $\G_{m,R}$.
and let $W$ be a function of weight $d$ on $V$. Also, we fix a finite subgroup $G\sub \GL(V)$ such that
$W$ is $G$-invariant, $G$ commutes with $\G_{m,R}$ and $G$ contains a fixed element $J\in \G_{m,R}$
of order $d$.

We define $\Ga\sub \GL(V)$ to be the algebraic subgroup generated by $G$ and by $\G_{m,R}$. 
There is a canonical exact sequence
$$1\to G\to \Ga\rTo{\chi} \G_m\to 1,$$
where $\chi$ restricts to the subgroup $\G_{m,R}$ as $\la\mapsto \la^d$.

As in \cite{PV}, we consider the moduli space of $\Ga$-spin structures: it classifies stable orbicurves $(C,p_1,\ldots,p_n)$ 
equipped with $\Ga$-principal bundle $P$ (our convention is that we have a right action of $\Ga$  on $P$), 
together with an isomorphism $\chi_*P\rTo{\sim} \om^{\log}_C\setminus 0$.
We can think of the latter isomorphism as a morphism $\chi_P:P\to \om^{\log}_C\setminus 0$ satisfying
$$\chi_P(x\ga)=\chi(\ga)\cdot \chi_P(x)$$
for $\ga\in \Ga$.

In addition to requiring the coarse moduli of $C$ to be Deligne-Mumford stable, we require that
for each marked point $p_i$ the morphism $B\Aut(p_i)\to B\Ga$ induced by $P$ is representable.
By looking at the corresponding embedding $\Aut(p_i)\simeq \Z/m_i\to \Ga$ defined up to a conjugacy, we get
a conjugacy class $\ga_i$ in $\Ga$. Thus, we get a decomposition of our moduli stack into a disjoint union
of open and closed substacks $\SS_g(\ga_1,\ldots,\ga_n)$.
As in \cite[Sec.\ 2.2]{FJR}, one shows that these are smooth and proper DM stacks with projective coarse moduli.

Let $\pi:\CC\to \SS_g(\ga_1,\ldots,\ga_i)$ be the universal curve over $\SS_g(\ga_1,\ldots,\ga_n)$, and let 
$\VV=\PP\times_\Ga V$ be the vector bundle over 
$\CC$ associated with the universal $\Ga$-spin structure $\PP$ via the embedding $\Ga\sub\GL(V)$.
Note that $\VV$ is equipped with a $\G_{m,R}$-action (through its action on $V$).

As in \cite{PV}, we also consider a Galois covering 
$\SS^{\rig}_g(\ga_1,\ldots,\ga_n)\to \SS_g(\ga_1,\ldots,\ga_n)$ corresponding to choices of a rigidification at every marked point. 
A {\it rigidification} is an isomorphism of the restriction of $P$ to 
$p_i/\Aut(p_i)\simeq B\lan \ga_i\ran$ with $V/\lan \ga_i\ran$
(viewed as a bundle over $B\lan \ga_i\ran$.
There is a natural simply transitive action of the group $\prod_i C_G(\ga_i)$ 
on the set of rigidifications at $p_1,\ldots,p_n$, where $C_G(\ga)\sub G$ is the centralizer of
$\ga\in G$.

\subsection{Construction}

Let us set for now $\SS=\SS^{\rig}_g(\ga_1,\ldots,\ga_n)$ and consider the pull-back of all the objects to $\SS$
(denoting them by the same symbols).

Note that we have a natural projection $V/\lan \ga_i\ran \to V^{\ga_i}$.
Thus, from rigidification structures we get morphisms
\begin{equation}\label{rig-i-map}
Z_i:p_i^*\VV\to V^{\ga_i}\ot \OO_{\SS}.
\end{equation}
Hence, by adjunction we can extend $\VV$ to an object 
$$[\VV,\Si]:=(\VV,(V^{\ga_i}\ot \OO_\SS),(Z_i))$$ 
of $\Coh(\CC,\Si)$.

On the other hand, we can combine $\chi_P$ with $W$ into a polynomial morphism
$$W_\VV:\VV=\PP\times_\Ga V\to \om^{\log}_{\CC/\SS}: (x,v)\mapsto W(v)\cdot \chi_P(x).$$
We can view it as a linear morphism of vector bundles on $\CC$,
$$W_\VV:S^\bullet(\VV)_d\to \om^{\log}_{\CC/\SS},$$
where we grade the symmetric algebra of $\VV$ using the $\G_{m,R}$-action on $\VV$.
Furthermore, this morphism is compatible with the morphisms \eqref{rig-i-map}, so that the following diagram is commutative
\begin{diagram}
p_i^* S^\bullet(\VV)_d&\rTo{p_i^*W_\VV}& p_i^*\om^{\log}_{\CC/\SS}\\
\dTo{S^\bullet(Z_i)}&&\dTo{}\\
S^\bullet(V^{\ga_i})_d\ot\OO_\SS&\rTo{W_i}& \OO_\SS
\end{diagram}
where $W_i=W|_{V^{\ga_i}}$.
This means that we have a morphism
\begin{equation}
(W_\VV,(W_i)):S^\bullet[\VV,\Si]_d\to [\om_{\CC/\SS}^{\log},\Si]
\end{equation}
in the category $\Qcoh(\CC,\Si)$ (where again we take the part of weight $d$ with respect to $\G_{m,R}$).
Next, we can take the derived push-forward with respect to the morphism of pairs
\eqref{projection-poset-eq}.
Together with \eqref{push-forward-tensor-eq} this gives us a morphism
\begin{equation}\label{S-V-mor}
S^\bullet(R\pi_*[\VV,\Si])_d\to R\pi_*S^\bullet[\VV,\Si]_d\to R\pi_*[\om_{\CC/\SS}^{\log},\Si]
\end{equation}
in $D\Qcoh(\SS,\SS)$.

Now let us set 
$$E:=Rj^!S^\bullet(R\pi_*[\VV,\Si])_d.$$
Applying $Rj^!$ to morphism \eqref{S-V-mor}, we obtain a morphism
$$E=Rj^!S^\bullet(R\pi_*[\VV,\Si])_d\to Rj^!R\pi_*[\om_{\CC/\SS}^{\log},\Si]\simeq R\pi_*\om_{\CC/\SS},$$
where the last isomorphism is \eqref{j!-om-log-eq}.
It is easy to see that it fits into a morphism of exact triangles
\begin{equation}\label{E-diagram}
\begin{diagram}
E&\rTo{}& S^\bullet(R\pi_*(\VV))_d &\rTo{}& \oplus_{i=1}^r S^\bullet(V^{\ga_i})_d\ot \OO_\SS\\
\dTo{}&&\dTo{}&&\dTo{(W_i)}\\
R\pi_*(\om_{\CC/\SS})&\rTo{}&R\pi_*(\om^{\log}_{\CC/\SS})&\rTo{}& \oplus_{i=1}^r\OO_\SS
\end{diagram}
\end{equation}
Combining it with the morphism of triangles \eqref{diff-res-diagram}, we get a commutative diagram
with the exact triangle in the first row
\begin{diagram}
S^\bullet(R\pi_*(\VV))_d &\rTo{}& \oplus_{i=1}^r S^\bullet(V^{\ga_i})_d\ot \OO_\SS&\rTo{}&E[1]\\
&& \dTo{\sum W_i}&&\dTo{\tau}\\
 &&\OO_\SS&\rTo{\id}&\OO_\SS
\end{diagram}
Dualizing we get a commutative diagram
\begin{diagram}
E^\vee[-1]&\rTo{}&\oplus_{i=1}^r S^\bullet(V^{\ga_i})_d^\vee\ot \OO_\SS&\rTo{}&
S^\bullet(R\pi_*(\VV))_d^\vee  \\
\uTo{\tau^\vee}&& \uTo{\sum W_i}\\
\OO_\SS&\rTo{\id}&\OO_\SS
\end{diagram}
This implies that the pull-back $Z^*(\bigoplus_i W_i)$ with respect to the morphism
\begin{equation}\label{Z-Rpi*-eq}
Z:[R\pi_*(\VV)]\to \prod_i V^{\ga_i}
\end{equation}
induced by \eqref{rig-i-map}, becomes zero in cohomology of the structure sheaf on $[R\pi_*(\VV)]$. 

In fact, we can realize this functon by an explicit coboundary. For this we need a realization of 
the above diagram in the homotopy category of complexes. As in \cite[Sec.\ 4.2]{PV}, the starting point is that
$R\pi_*(\VV)$ can realized ($\G_{m,R}$-equivariantly) by a complex of the form $[A\to B]$ in such a way
that the morphism \eqref{Z-Rpi*-eq} is realized by a surjective morphism 
$A\to \bigoplus_{i=1}^r V^{\ga_i}\ot\OO_S$. Then 
the first line of the diagram \eqref{E-diagram} 
can be realized by a short exact sequence of complexes
$$0\to \ker(S^\bullet(Z)_d) \to S^\bullet(A\to B)_d \rTo{S^\bullet(Z)_d} \bigoplus_{i=1}^r S^\bullet(V^{\ga_i})_d\ot \OO_S\to 0$$
where the complex $S^\bullet(A\to B)_d$, concentrated in degrees $[0,\rk(B)]$, has form
$$S^\bullet(A)_d\to (S^\bullet(A)\ot B)_d\to (S^\bullet(A)\ot {\wedge}^2 B)_d\to \ldots
$$
Using this we get a canonical quasi-isomorphism of $E$ with the bounded complex of vector bundles
$$K^\bullet:=\Cone(S^\bullet(R\pi_*(\VV))_d \to \oplus_{i=1}^r S^\bullet(V^{\ga_i})_d\ot \OO_\SS)[-1].$$

Now we want to realize the morphism $\tau:E\to \OO_S[-1]$ in the derived category
by a morphism $K^\bullet\to \OO_S[-1]$ in the homotopy category of complexes. 

By changing $[A\to B]$ to a quasi-isomorphic complex $[\ov{A}\to \ov{B}]$ one
can achieve that for $i\ge 1$ the terms $K^i$ satisfy $\Ext^{>0}(K^i,\OO_S)=0$ (see \cite[Lem.\ 4.2.5]{PV}). 
This implies that morphisms $K\to \OO_S[-1]$
in the homotopy category of complexes and in the derived category are the same.

The dual of this morphism can be interpreted as a canonical homotopy (up to a homotopy between homotopies) $f_{-1}$ between
the function $Z^*(\bigoplus_i W_i)$ on $[R\pi_*\VV]$ and $0$. As we have seen in Example \ref{cob-ex}, this corresponds
to a structure $\de=d-f_{-1}\cdot \id$ of a dg-matrix factorization of $-Z^*(\bigoplus_i W_i)$ on the structure sheaf of
$[R\pi_*\VV]$. 

Furthermore, it carries an equivariant structure with respect to the action of the center $Z(\Ga)$ of $\Ga$
(acting trivially on the base) and with respect to $\prod_i C_G(\ga_i)$ (changing the rigidifications).

\subsection{Properties}

The first important property is that our dg-matrix factorization over $[R\pi_*\VV]$ is supported on the zero section in $[R\pi_*\VV]$.
Since each $W_i$ is non-degenerate, we know that the support belongs to the zero locus of $Z^*(\bigoplus_i W_i)$.
Thus, we reduce to considering the following situation. Let $C$ be a curve, $\VV$ be a vector bundle over $C$,
equipped with a $\G_m$-equivariant structure (where $\G_m$ acts trivially on $C$).
Assume also we have a polynomial morphism $W_\VV:\VV\to \om_C$, homogeneous of degree $d$, such that
over an open dense subset of $C$ there exists a trivialization $\VV\simeq V\ot\OO_C$ (compatible with the
$\G_m$-action) such that $W_\VV$ is induced by our polynomial $W$ on $V$. Then we have the induced
polynomial function of degree $-1$ on the dg-affine space $[H^0(C,\VV)\oplus H^1(C,\VV)[-1]]$, induced
by $W_\VV$ and by the identification $H^1(C,\om_C)\simeq\C$. We claim that it is supported at the origin.
Indeed, we start by observing that the preimage of the origin under the gradient morphism 
$\Delta W: V\to V^\vee$ is still the origin (since $W$ is non-degenerate). 
From this we get the similar assertion about the preimage of the zero section
under the relative gradient morphism $\Delta W_\VV: \VV\to \VV^\vee\ot\om_C$. Finally, we note that
the support of our function on $[H^0(C,\VV)\oplus H^1(C,\VV)[-1]]$ coincides with the vanishing locus of
the polynomial morphism
$$H^0(C,\VV)\to H^0(\VV^{\vee}\ot\om_C)\simeq H^1(C,\VV)^\vee$$
induced by the relative gradient map. This implies our claim.

Next, the key gluing property satisfied by the fundamental matrix factorizations (cf. \cite[Sec.\ 5.2, 5.3]{PV}) holds
in the situation when we consider two natural families of orbicurves $\wt{C}\rTo{\wt{\pi}} S$, $C\rTo{\pi} S$, over
$$S:=S^{\rig}_{g_1}(\ga_1,\ldots,\ga_{n_1},\ga)\times S^{\rig}_{g_2}(\ga'_1,\ldots,\ga'_{n_2},\ga^{-1}),$$
where $\wt{C}$ is the disconnected curve and $C$ is obtained by gluing two points into a node
(there is also a similar picture corresponding to a non-disconnecting node).
We denote by $f:\wt{C}\to C$ the gluing morphism.

In this setting there are natural $\Ga$-spin structures $\wt{P}$ (resp., $P$) over $\wt{C}$ (resp., $C$),
where $P$ is obtained by gluing fibers of $\wt{P}$ over the two points that are glued into a node, using
the rigidifications and the square root of $J$, $J^{1/2}\in \G_{m,R}$ such that $\chi(J^{1/2})=-1$ 
(see \cite[Sec.\ 5.2]{PV}).
The main compatibility between the push-forwards of the corresponding vector bundles $\wt{\VV}$ and $\VV$ is 
given by the cartesian diagram
\begin{diagram}
[R\pi_*\VV] &\rTo{}& V^{\ga}\\
\dTo{}&&\dTo{\De^{J^{1/2}}}\\
[R\wt{\pi}_*\wt{\VV}] &\rTo{}& V^{\ga}\times V^{\ga^{-1}}
\end{diagram}
where $\De^{J^{1/2}}:V^{\ga}\to V^{\ga}\times V^{\ga^{-1}}$ is the twisted diagonal map: $x\mapsto (x,J^{1/2}x)$.
Furthermore, the natural dg-matrix factorization on $[R\pi_*\VV]$ is identified with the pull-back of the one on 
$[R\wt{\pi}_*\wt{\VV}]$.


Recall that in \cite{PV}, in the case when $G$ is contained in an algebraic torus acting on $V$, we used the fundamental matrix factorizations to construct cohomological field theories
associated with $(W,G)$ by viewing them as kernels for Fourier-Mukai functors and passing to Hochschild homology. 
Similarly, in the case of nonabelian $G$ one can use
the above construction to get a cohomological field theory (with coefficients in $\C$)
and show some of the properties predicted in \cite{FJR2}. We will study the resulting cohomological field theories elsewhere.



\begin{thebibliography}{9}
\bibitem{CF-K} I.~Ciocan-Fontanine, M.~Kapranov, {\it Virtual fundamental classes via dg-manifolds}, 
Geom. Topol. 13 (2009), 1779--1804.
\bibitem{EP} A.~I.~Efimov, L.~Positselski,
{\it Coherent analogues of matrix factorizations and relative singularity categories}, Algebra Number Theory 9 (2015), 1159--1292. 
\bibitem{FJR} H.~Fan, T.~Jarvis, Y.~Ruan, {\it The Witten equation, mirror symmetry, and quantum singularity theory}, 
Ann. of Math. (2) 178 (2013), 1--106.
\bibitem{FJR2} H.~Fan, T.~Jarvis, Y.~Ruan,
{\it A mathematical theory of the gauged linear sigma model}, Geom. Topol. 22 (2018), 235--303.
\bibitem{lunts-poset} V.~Lunts,
{\it Categorical resolutions, poset schemes, and Du Bois singularities}, IMRN 19 (2012), 4372--4420.
\bibitem{PV} A.~Polishchuk, A.~Vaintrob, {\it Matrix factorizations and cohomological field theories}, 
J. Reine Angew. Math. 714 (2016), 1--122.
\end{thebibliography}
\end{document}